\documentclass{article}
\usepackage{maa-monthly}

\newcommand\Defn[1]{\textbf{\color{black}#1}}
\newcommand\Def[1]{\Defn{#1}}
\newcommand\R{\mathbb{R}}
\newcommand\Rnn{\R_{\ge0}}
\newcommand\Z{\mathbb{Z}}
\newcommand\Znn{\mathbb{Z}_{\ge0}}
\newcommand\G{\mathcal{G}}%
\newcommand\B{\mathcal{B}}%
\newcommand\J{\mathbf{J}}
\newcommand\Inner[1]{\langle{#1}\rangle}%
\newcommand\1{\mathbf{1}}
\renewcommand\emptyset{\varnothing}

\DeclareMathOperator{\conv}{conv}

\theoremstyle{theorem}
\newtheorem{theorem}{Theorem}
\newtheorem{lemma}{Lemma}
\newtheorem{proposition}{Proposition}

\theoremstyle{definition}

\begin{document}

\title{The Martin Gardner Polytopes}
\markright{The Martin Gardner Polytopes}
\author{Kristin Fritsch, Janin Heuer, Raman Sanyal, and Nicole Schulz}

\maketitle

\begin{abstract}
    In the chapter ``Magic with a Matrix'' in \emph{Hexaflexagons and Other
    Mathematical Diversions} (1988), Martin Gardner  describes a delightful
    ``party trick'' to fill the squares of a $d$-by-$d$ chessboard with
    nonnegative integers such that the sum of the numbers covered by any
    placement of $d$ nonthreatening rooks is a given number $N$.
    We consider such chessboards from a geometric perspective
    which gives rise to a family of lattice polytopes. The polyhedral
    structure of these Gardner polytopes explains the underlying trick and
    enables us to count such chessboards for given $N$ in three different
    ways. We also observe a curious duality that relates Gardner polytopes to
    Birkhoff polytopes.
\end{abstract}

\section{Introduction.}\label{sec:intro}

In \emph{Magic with a Matrix}~\cite[Chapter~2]{Gardner} Martin Gardner describes
the following delightful ``party trick.'' You ask the audience for two
positive integers $d$ and $N$. After a moment's thought, you fill the 
squares of a $d$-by-$d$ chessboard with nonnegative integers such
that the sum of the numbers covered by any placement of $d$ nonthreatening
rooks is $N$. For example, given $d=5$ and $N=57$ you might produce the
following chessboard (that we took from~\cite{Gardner}):
\[
    \begin{array}{ccccc}
        19 & 8  & 11 & 25 & 7 \\
        12 & 1  & 4  & 18 & 0 \\
        16 & 5  & 8  & 22 & 4 \\
        21 & 10 & 13 & 27 & 9 \\
        14 & 3  & 6  & 20 & 2
    \end{array}
\]
The trick to pick the numbers is so unsettlingly simple that with some
training, you will be able to fill the chessboard in no time. The audience is
left to wonder how you could do that so fast and the more mathematically
inclined members of the audience will automatically start asking themselves
how such fillings can be constructed systematically. We encourage readers
to do so for themselves before continuing.

In this note we want to approach this question from a geometric angle. This
reveals beautiful connections to gems of the theory of lattice polytopes. As a
first step, note that we might as well fill a $d$-by-$d$ table with
nonnegative \emph{real} numbers with the same property. We therefore call a
matrix $A \in \R^{d \times d}$ a \Defn{$\boldsymbol G$-matrix} of value $N \in \R_{\ge0}$
if all entries are nonnegative and
\begin{equation}\label{eqn:G}
    A_{1,\sigma(1)} + A_{2,\sigma(2)} + \cdots + A_{d,\sigma(d)}  \ = \ N
\end{equation}
for every permutation $\sigma$ of $[d] := \{1,\dots,d\}$. Note that the zero
matrix is the only $G$-matrix of value $N = 0$. If we let  $\G_d
\subseteq \R^{d \times d}$ be the collection of $G$-matrices of value $1$,
then $A$ is a $G$-matrix of value $N > 0$ if and only if $\frac{1}{N}  A \in
\G_d$. The set $\G_d$ is given by linear equations and inequalities and is
thus a convex polyhedron. To see that $\G_d$ is in fact a polytope, that is, a
bounded polyhedron, we simply note that the entries of $G$-matrices of value
$1$ inevitably are bounded from above by $1$ and hence $\G_d$ is contained
in the unit cube $[0,1]^{d\times d} \subseteq \R^{d\times d}$. For obvious
reasons, we call the polytopes $\G_d$ the \Def{Gardner polytopes}.

In the next section, we will study the dimension, the vertices, and the
general polyhedral structure of $\G_d$. This basically addresses the question
of how to construct \mbox{$G$-matrices.} In Section~\ref{sec:count}, we want
to count $G$-matrices with integer entries for a fixed value $N \in
\Z_{\ge0}$. That is, we want to count the points in $N \cdot \G_d \cap \Z^{d
\times d}$ and we give three representations of that counting function as a
polynomial in $N$ of degree $2d-2$. In the last section we observe a curious
duality that relates the Gardner polytope $\G_d$ to the Birkhoff polytopes
$\B_d$, the polytope of doubly stochastic $d$-by-$d$ matrices, and that
deserves more study. For more on polytopes and their combinatorics we refer
to~\cite{ziegler}.

After the first revision of our paper, Professor Richard Stanley kindly informed
us that he too considered such matrices (which he called \emph{antimagic
squares}) and an account of his findings can be found  in his marvelous
book~\cite[Exercise 4.53]{stanley}.  While the results (and even the notation)
are quite similar, our approach is firmly rooted in geometry (polytopes and
subdivisions), which gives new insights into
$G$-matrices/antimagic squares.

\section{Gardner polytopes.}\label{sec:polytope}
According to the Minkowski--Weyl theorem~\cite[Theorem~1.1]{ziegler}, if
$\G_d$ is a polytope, then $\G_d$ is the intersection of all convex sets
containing a finite set of points $S \subset \R^{d \times d}$. We write this
as $\G_d = \conv(S)$ and call $\G_d$ the \Def{convex hull} of $S$. There is a
unique inclusion-minimal set $S$ with $\G_d = \conv(S)$ and its elements are
called the \Def{vertices} of $\G_d$. Determining the vertices of $\G_d$ is
what we do first.

For $i = 1,\dots,d$, let $R_i, C_i \in \R^{d\times d}$ be the matrices whose
entries in the $i$th row, respectively, $i$th column are $1$ and the
remaining entries are all $0$. Any placement of $d$ nonthreatening rooks has
exactly one rook in the $i$th row and column, and thus $R_1,\dots,R_d,
C_1,\dots,C_d \in \G_d$. 

\begin{theorem}\label{thm:vertices}
    The vertices of $\G_d$ are exactly the matrices $R_1,\dots,R_d,
    C_1,\dots,C_d$.
\end{theorem}

To prove this result, let $A \in \G_d$ be arbitrary. For $j = 1,2,\dots,d$, let
$\lambda_j$ be the minimal entry in the $j$th column of $A$. Then $A' := A -
\lambda_1 C_1 - \lambda_2 C_2 - \cdots - \lambda_d C_d$ is a nonnegative
matrix and, in fact, a $G$-matrix of value $1 - \lambda_1 - \cdots -
\lambda_d$.  Analogously, let $\mu_i$ be the minimum in the $i$th row of $A'$
for $i=1,2,\dots,d$ and set $A'' := A' - \sum_i \mu_i R_i$, which is a
nonnegative $G$-matrix by the same reasoning. We claim that $A''$ is the zero
$G$-matrix.  Assume on the contrary that $A''_{ij} > 0$ for some $i,j \in
[d]$. By construction there are $k,l \in [d]$ such that $A''_{il} = A''_{kj} =
0$.  Choosing the permutations $\sigma,\sigma'$ with $(\sigma(i),\sigma(k)) =
(j,l)$ and $(\sigma'(i),\sigma'(k)) = (l,j)$ and $\sigma(h) = \sigma'(h)$ for
$h \in [d] \setminus \{i,k\}$, we infer from~\eqref{eqn:G} that 
\[
    0 \ = \ A''_{il} + A''_{kj} \ = \ A''_{ij} + A''_{kl} \, ,
\]
but $A''_{ij} > 0$ then implies $A''_{kl}< 0$, which is a contradiction. Thus
\begin{equation}\label{eqn:rep}
    A \ = \ 
    \lambda_1 C_1 + \cdots + \lambda_d C_d + 
    \mu_1 R_1 + \cdots + \mu_d R_d  \, ,
\end{equation}
and $\lambda_i, \mu_j \ge 0$ for all $i,j$ and $\lambda_1  + \cdots + \lambda_d  + \mu_1  + \cdots + \mu_d  = 1$.
This argument reveals the secret of $G$-matrices: every
$G$-matrix is simply an addition table for the column labels $\lambda_j$ and
the row labels $\mu_i$. Here is the situation for the example in the
Introduction:
\[
    \begin{array}{c|ccccc}
     + &   12 & 1  & 4  & 18 & 0 \\
     \hline
     7 &   19 & 8  & 11 & 25 & 7 \\
     0 &   12 & 1  & 4  & 18 & 0 \\
     4 &   16 & 5  & 8  & 22 & 4 \\
     9 &   21 & 10 & 13 & 27 & 9 \\
     2 &   14 & 3  & 6  & 20 & 2
    \end{array}
\]
The \Def{affine hull} of $\G_d$, the inclusion-minimal affine subspace
containing $\G_d$, is given by
\begin{equation}\label{eqn:aff}
    \begin{aligned}
        A_{11} + A_{12} + \cdots + A_{dd} \ = \ d  & \quad \text{and} \\
    A_{ij} + A_{kl} \ = \ A_{il} + A_{kj} & \quad \text{for }
    i \neq k \text{ and } j \neq l \, .
    \end{aligned}
\end{equation}
In the language of~\cite[Section~5.3]{MS}, the latter set of equations states
that $A$ is a \emph{tropical} rank-$1$ matrix.

We still have to argue that we cannot forgo any of the proposed matrices and
hence $R_1,\dots,R_d,C_1,\dots,C_d$ are indeed the vertices of $\G_d$. 
We claim that 
\begin{equation}\label{eqn:circuit}
    R_1 + \cdots + R_d \ = \ \J \ = \ C_1 + \cdots + C_d \, ,
\end{equation}
where $\J$ is the \Def{all-ones matrix}, gives the unique affine dependence
among these matrices up to scaling. Let $\alpha_i, \beta_j \in \R$ not all
zero be such that 
\[
    \alpha_1 R_1 + \cdots + \alpha_d R_d 
    \ = \
    \beta_1 C_1 + \cdots + \beta_d C_d \, .
\]
By inspecting the entry $(i,j)$, we see that $\alpha_i = \beta_j$ and hence we
can choose $\alpha_i=\beta_j = 1$ for all $i,j$. Now if, say, $R_1$ is
not a vertex of $\G_d$, then $R_1$ is a convex linear combination of
$R_2,\dots,R_d,C_1,\dots,C_d$, which would yield an affine dependence other
than~\eqref{eqn:circuit}.

Equation~\eqref{eqn:circuit} marks the set $\{R_1,\dots,R_d,C_1,\dots,C_d\}$
as a \Def{circuit}, i.e.,\ an inclusion-minimal collection of affinely
dependent points. This makes $\G_d$ quite simple from a 
polytope point of view. Let us write $\Inner{A,B} = \mathrm{tr}(A^tB) =
\sum_{i,j} A_{ij}B_{ij}$ for the standard inner product on $\R^{d \times d}$.
Then
\begin{align*}
    P_R \ &:= \ \conv ( R_1,\dots, R_d ) \ = \ \{ A \in \G_d : \Inner{A,C_j}
    \ = \ 1 \text{ for } j = 1,\dots,d \} \ \text{ and }\\
    P_C \ &:= \ \conv ( C_1,\dots, C_d ) \ = \ \{ A \in \G_d : \Inner{A,R_i}
    \ = \ 1 \text{ for } i = 1,\dots,d \} \, .
\end{align*}
$P_R$ as well as $P_C$ are given as the convex hull of $d$ affinely
independent points and hence are \Def{simplices} of dimension $d-1$.
 Moreover~\eqref{eqn:circuit} yields that $P_R
\cap P_C = \{ \frac{1}{d} \J \}$. The polytope $\G_d$ is the convex hull of the
union of $P_R$ and $P_C$, which is called the \Def{direct sum} (or sometimes
the \Def{free sum}) of $P_R$ and $P_C$. In summary, $\G_d$ is a
$(2d-2)$-dimensional polytope on $2d$ vertices and, relative to the affine
hull~\eqref{eqn:aff}, the $d^2$ inequalities $A_{ij} \ge 0$ give an
irredundant description.

\section{Counting $G$-matrices.}\label{sec:count}
We want to count the number of integer \mbox{$G$-matrices} of fixed size $d$ and
value $N \ge 0$; let $g_d(N)$ be this number. In pursuing our geometric
approach to $G$-matrices, we make the following observation.

\begin{lemma}
    For $N \in \Znn$, $g_d(N) = |N\cdot \G_d \cap \Z^{d \times d}|$.
\end{lemma}

Why should that make determining $g_d(N)$ easier?  Counting lattice points in
polytopes has a long (and ongoing) history, in particular if the lattice
points in the polytope under scrutiny have combinatorial meaning;
see~\cite{crt} for much more on this.  In 1962 Eug\`ene
Ehrhart~\cite{ehrhartpolynomial} discovered that counting lattice points in
dilates of polytopes is particularly nice if the underlying polytope has all
its vertices in the integer lattice; we therefore call this a \Defn{lattice
polytope}.

\renewcommand\P{\mathcal{P}}%
\begin{theorem}[Ehrhart]\label{thm:ehrhart}
    Let $\P \subset \R^m$ be a lattice polytope of dimension $D$. Then the
    function $E_\P(n) := |n\P \cap \Z^m|$ agrees with a polynomial of degree
    $D$---the
    \Def{Ehrhart polynomial} of $\P$---for all $n \in \Z_{\ge0}$.
\end{theorem}

Together with Theorem~\ref{thm:vertices} we readily conclude that $g_d(N)$
is a polynomial in $N$ of degree $2d-2$.

There are many perspectives on and proofs of Theorem~\ref{thm:ehrhart};
see~\cite{BR,crt} for a selection. A sound approach is by way of subdivisions
of polytopes, and we too will subdivide $\G_d$ into simpler polytopes to give
\emph{three} different representations of $g_d(N)$ given in
equations~\eqref{eqn:1st}, \eqref{eqn:2nd}, and \eqref{eqn:3rd}.

Consider the polytopes
\[
    P_k \ := \ \conv(R_1,\dots,R_{k-1},R_{k+1},\dots,R_d,C_1,\dots,C_d)
\]
for $k = 1,\dots,d$. Since the vertices of $\G_d$ form a circuit, the vertices
of $P_k$ are affinely independent and hence $P_k$ is a $(2d-2)$-dimensional
simplex for every $k$. For a polytope $\P$, a \Def{face} of $\P$ is a subset
$F \subseteq \P$ such that for $p,q \in \P$ the midpoint $\frac{p+q}{2}$ is
contained in $F$ only if $p,q \in F$. If $\P$ is a simplex, it is not
difficult to check that faces are precisely convex hulls of arbitrary subsets
of the vertices of $\P$.

\begin{proposition}\label{prop:triang}
    The simplices $P_1,\dots,P_d$ form a \Def{triangulation} of $\G_d$, that
    is, 
    \[
        \G_d \ = \  P_1 \cup P_2 \cup \cdots \cup P_d
    \] 
   and $P_i \cap P_j$ is a face of both $P_i$ and $P_j$ for all $1 \le i < j
   \le d$.
\end{proposition}

The first condition follows from our representation~\eqref{eqn:rep}. Indeed,
by choosing the column minima first, at least one of the rows of the matrix
$A'$ has a zero entry and thus not all the $\mu_i$ are positive. But if $\mu_k
= 0$ in~\eqref{eqn:rep}, then $A \in P_k$. To verify the second claim, we note
that if a point $p \in P_i \cap P_j$ has different representations with
respect to the vertices of $P_i$ and $P_j$, then this gives an affine
dependence different from~\eqref{eqn:circuit}, which contradicts its
uniqueness. Hence 
\begin{equation}\label{eqn:intersect}
    P_i \cap P_j \ = \ \conv(C_1,\dots,C_d, R_k : k \in [d] \setminus \{i,j\})
\end{equation}
is a face of both polytopes.

Of course, omitting $C_i$ one at a time from the convex hull yields another
triangulation and, by virtue of circuits, these are the only two
triangulations of $\G_d$ without new vertices; see, for
example,~\cite[Section~2.4]{Triang}.

To understand the arithmetic of, say, $P_1$ better, we make use of a linear projection
$\pi : \R^{d \times d} \to \R^{2d-2}$. If $A \in \R^{d \times d}$ is a matrix
with first row $(r_1,r_2,\dots,r_d)$ and first column $(c_1 = r_1,c_2,\dots,c_d)$,
then we set
\[
    \pi(A) \ := \ (r_2,r_3,\dots,r_{d}, c_2 - c_1, c_3 - c_1, \dots, c_d -
    c_1) \, .
\]
Let $e_i$ be the $i$th standard vector in $\R^{2d-2}$ and set $e_0 := 0$.
Then under this projection, we have $\pi(C_j) = e_{j-1}$ for $j \ge 1$ and
$\pi(R_{i+1}) = e_{d-1+i}$ for $i =1,\dots,d-1$. Thus $\pi(P_1) =
\Delta_{2d-2} := \conv(0,e_1,\dots,e_{2d-2})$ is the \Def{standard simplex} in
$\R^{2d-2}$. The projection $\pi$ gives a linear isomorphism between $P_1$ and
$\Delta_{2d-2}$, and for $N \in \Znn$ a point $A \in N\cdot P_1$ is a lattice
point if and only if $\pi(A) \in N\cdot\Delta_{2d-2}$ is. This shows that 
\[
    E_{P_i}(N) \ = \ E_{\Delta_{2d-2}}(N)
\]
for all $i = 1,\dots,d$ and $N \in \Znn$.  Simplices with this Ehrhart
polynomial are called \Def{unimodular} and are characterized by having the
minimal volume among all simplices of the same dimension with vertices in the
integer lattice. The Ehrhart polynomial of $\Delta_{m-1}$ is
$E_{\Delta_{m-1}}(n) = \binom{n+m-1}{m-1}$ (go on, try it!).  At this point,
we can appeal to the principle of inclusion-exclusion for the triangulation of
Proposition~\ref{prop:triang}. For $K \subseteq [d]$ note that $\bigcap_{i \in
K} P_i \ = \ \conv(C_1,\dots,C_d, R_i : i \not\in K) \cong \Delta_{2d -
|K|-1}$. Hence we get our first expression for $g_d(N)$:
\begin{equation}\label{eqn:1st}\tag{\textbf{1st}}
    g_d(N) \ = \ \sum_{k = 1}^d (-1)^{k-1} \binom{d}{k} \binom{N + 2d - k -1
    }{2d - k-1} \, .
\end{equation}

Proposition~\ref{prop:triang} actually states that $\G_d$ has a
\Def{unimodular triangulation}, i.e., a triangulation into unimodular
simplices. A unimodular triangulation is a luxurious thing to have as the
Ehrhart polynomial depends only on the combinatorics of the triangulation
(see~\cite[Section~5.5]{crt}) which in our case is pretty simple.

To see how this works, we define the \Def{relative interior} $\P^\circ$ of a
polytope $\P \subset \R^m$ as the interior of $\P$ relative to its affine
hull. It is a basic fact that every polytope is the disjoint union of the
relative interiors of its faces. We want to use this together with
Proposition~\ref{prop:triang} to write $\G_d$ as the disjoint union of the
faces of the $P_i$. The obvious benefit is that we can avoid
inclusion-exclusion by counting lattice points in the various relatively open
pieces.  Since the $P_i$ meet in faces, we have to determine which subsets of
$\{R_1,\dots,R_d,C_1,\dots,C_d\}$ form a face of one of the $P_i$. For $I,J
\subseteq [d]$, we write $P_{I,J}$ for the convex hull of the points $\{R_i :
i \in I\}$ and $\{C_j : j \in J \}$. Then $P_{I,J}$ occurs as a face of one of
the $P_i$ if and only if $I \neq [d]$.  The simplex $P_{I,J}$ is lattice
isomorphic to $\Delta_{m-1}$ with $m = |I|+|J|$ and, either by appealing to
Ehrhart--Macdonald reciprocity~\cite[Section~5.4]{crt} or simply by writing it
out, we observe that $E_{\Delta^\circ_{m-1}}(n) = \binom{n-1}{m-1}$.  We
conclude that $|N \cdot P_{I,J}^\circ \cap \Z^{d \times d}| =
\binom{N-1}{m-1}$ for all $N \in \Z_{\ge 0}$ and counting the number of such
sets $I,J$ gives the second expression of $g_d(N)$:
\begin{equation} \label{eqn:2nd}\tag{\textbf{2nd}}
    g_d(N) \ = \ \sum_{m = 1}^{2d-1} \left[\binom{2d}{m} -
    \binom{d}{m-d}\right] \binom{N-1}{m-1} \, .
\end{equation}
Our reasoning gives a first idea of the interplay of geometry and
combinatorics: for any $D$-dimensional lattice polytope $\P$ we may write its
Ehrhart polynomial in the form $E_\P(n) = \sum_{i=1}^{D+1} f_{i-1}^\ast(\P)
\binom{n-1}{i-1}$. Then any unimodular triangulation of $\P$, in case one
exists, has to have $f_{i-1}^\ast$ many simplices of dimension $i-1$; for more
see~\cite{Breuer}.

A more economical way to avoid inclusion-exclusion is by using pieces that are
neither closed nor open.  Let $\P = \conv(v_1,\dots,v_{m})$ be a
$(m-1)$-dimensional simplex.  Every point $q \in \P$ is of the form 
\begin{equation}\label{eqn:conv}
    q \ = \ \lambda_1 v_1 +  \cdots + \lambda_{m} v_{m}
\end{equation}
for some unique $\lambda_1,\dots,\lambda_{m} \ge 0$ with $\lambda_1
+  \cdots +\lambda_{m} = 1$. For $U \subseteq 
\{v_1,\dots,v_m\}$, let us write $H_U \P$ for the set of points $q \in \P$ for
which in~\eqref{eqn:conv} we have $\lambda_i > 0$ whenever $v_i
\in U$. We call $H_U\P$ a \Def{half-open} simplex. The two extreme cases are
$H_\emptyset \P = \P$ and $H_{\{v_1,\dots,v_m\}}\P = \P^\circ$. If $\P$ is a unimodular simplex, then one checks that 
\[
    E_{H_U\P}(n) \ = \ \binom{n-1+m - |U|}{m-1} \, .
\]
So, to get our third expression for $g_d(N)$, we want to find suitable
$U_1,\dots,U_d$ such that 
\[
    \G_d \ = \ 
    H_{U_1}P_1 \uplus
    H_{U_2}P_2 \uplus
    \cdots \uplus
    H_{U_d}P_d \, ,
\]
which we call a \Def{half-open decomposition}.  We claim that $U_i := \{
R_1,\dots,R_{i-1} \}$ for $i=1,\dots,d$ do the job.  That the union on the
right-hand side is disjoint follows directly from~\eqref{eqn:intersect}. So we
only need to argue that $\G_d$ is covered by the half-open simplices. But for
any $A \in \G_d$ we devised a canonical representation given
in~\eqref{eqn:rep}. Letting $k \ge 1$ be minimal with $\mu_k = 0$, we observe
that $A \in H_{U_k}P_k$ which finishes the argument.  Since $|U_i| = i-1$, we
obtain our third expression:
\begin{equation}\tag{\textbf{3rd}}\label{eqn:3rd}
    g_d(N)  =  \sum_{j=0}^{d-1} \binom{N + 2d - 2 - j}{2d-2}  =  
    \binom{N + 2d-1}{2d-1} - 
    \binom{N + d-1}{2d-1} \, .
\end{equation}
Half-open decompositions, introduced in~\cite{KV}, are a powerful tool for
proving structural results as well as for practical computations;
see~\cite[Section~5.3]{crt}. An interesting property of this last expression
of $g_d(N)$ is that all its zeros are either negative integers or have real
part $-\frac{d}{2}$, as can be inferred from~\cite[Theorem~3.2]{StanCyc}. So
the Ehrhart polynomials of the Gardner polytopes satisfy a \emph{Riemann
hypothesis} in the sense of~\cite{BCKV}; see also~\cite{HKM}.

\section{An interesting duality.}\label{sec:duality}
A matrix $B \in \R^{d \times d}$ is \Def{doubly stochastic} if $B$ is
nonnegative and all row and column sums are equal to $1$.  The collection of
doubly stochastic $d$-by-$d$ matrices is called the \Def{Birkhoff polytope}
$\B_d \subset \R^{d \times d}$. This is a convex polytope of dimension
$(d-1)^2$ and, according to the Birkhoff--von Neumann theorem, the vertices of
$\B_d$ are precisely the permutation matrices;
see~\cite[Section~II.5]{barvinok}. For a permutation $\sigma$ of $[d]$ the
corresponding permutation matrix $P_\sigma \in \R^{d \times d}$ is the
$0/1$-matrix with $(P_\sigma)_{ij} = 1$ if and only if $\sigma(i) = j$.  

The Gardner and Birkhoff polytopes satisfy an interesting duality:
\begin{align*}
\G_d \ &= \ \R^{d \times d}_{\ge0}  \cap  
\{ A : \Inner{A,P_\sigma}  = 1 \text{ for } \sigma \text{ permutation} \} \, , \\
\B_d \ &= \ \conv(P_\sigma : \sigma \text{ permutation}) \, ,
\intertext{
and, equivalently, }
\G_d \ &= \ \conv(R_i,C_j  :  i,j = 1,2,\dots,d ) \, , \\
\B_d \ &= \ \R^{d \times d}_{\ge0}  \cap  
\{ B : 
\Inner{B,R_i}  = 
\Inner{B,C_j}  = 1 \text{ for }  i,j = 1,2,\dots,d \} \, .
\end{align*}
\newcommand\Q{\mathcal{Q}}%
In this section we outline this duality of polytopes in general. We call $\P,\Q
\subset \R^D$ a \Def{Gale-dual} pair of polytopes if 
\begin{align*}
    \P & \ = \ \Rnn^D \ \cap \ 
    \{ x \in \R^D : \Inner{x,y} = 1 \text{ for } y \in V(\Q) \} \, , \\
    \Q & \ = \ \Rnn^D \ \cap \ 
    \{ y \in \R^D : \Inner{x,y} = 1 \text{ for } x \in V(\P) \}  \, , 
\end{align*}
where $V(\P)$ and $V(\Q)$ are the sets of vertices of $\P$ and $\Q$,
respectively. The naming comes from a certain reminiscence of ``Gale duality''
that we explain now.

For an affine subspace $L \subseteq \R^D$, let $L^\dag$ be the affine subspace
of all $y \in \R^D$ such that 
\[
    \Inner{x,y} \ = \ 1 \quad \text{ for all } \quad x \in L \, .
\]
Of course, $L^\dag \neq \varnothing$ if and only if $0 \not\in L$. 
If $L$ is an affine space not containing the origin, then $L = q + U$, where
$U \subseteq \R^D$ is a linear subspace and $q$ is a point with $U \subseteq
q^\perp$. It is easy to verify that 
\[
    L^\dag \ = \ \tfrac{1}{\|q\|^2} q + (U^\perp \cap q^\perp) \, .
\]
In particular, if we set $V = U^\perp \cap q^\perp$, then $U \oplus V =
q^\perp$ and $U \perp V$.  If $q = \1 := (1,\dots,1)$, then let the rows of $C
\in \R^{d \times D}$ and $G \in \R^{(D-d-1)\times D}$ be the bases for $U$ and
$V$, respectively. Then the columns of $C$ give a centered point configuration
with Gale transform $G$ and conversely; see~\cite[Section~5.6]{matousek}. 

In general, if $q >0$, then $\P = \Rnn^D \cap L$ and $\Q = \Rnn^D \cap L^\dag$
are a Gale-dual pair and plenty of Gale-dual pairs of polytopes can be
constructed from this simple recipe.  However, the Gale-dual pair $\G_d, \B_d$
is special. For starters, both $\G_d$ and $\B_d$ are lattice polytopes.  Both
polytopes are \Def{Gorenstein} of index $d$: for $\P = \G_d$ or $\P = \B_d$,
the all-ones matrix $\J$ is the unique lattice point in the interior of $d
\cdot \P$ and for $N \ge d$ a matrix $A$ is a lattice point in the interior of
$N \cdot \P$ if and only if $A - \J$ is a lattice point in $(N-d) \cdot \P$.
See~\cite[Section~6.E]{BrunsGubeladze} for much more on Gorenstein polytopes.
Moreover, the presentations of $\G_d$ and $\B_d$ above show that they are
\Def{compressed} polytopes in the sense of~\cite{sullivant}: $\B_d$ as well as
$\G_d$ are lattice polytopes given as the intersection of the unit cube with
an affine subspace.  Such polytopes have many desirable properties (e.g.,
every pulling triangulation is unimodular) but they are rare! It would be very
interesting to know if there are any other (nontrivial) examples/families of
Gale-dual pairs with all/any of these properties.

\begin{acknowledgment}{Acknowledgments.}
    We thank Professor Stanley for pointing out~~\cite[Exercise 4.53]{stanley}
    as well as the observation concerning the roots of $g_d(N)$ at the end of
    Section~\ref{sec:count}.

This paper grew out of a project of the course \emph{Polytopes,
Triangulations, and Applications} at Goethe-Universit\"{a}t 
Frankfurt in March 2018. We thank Sebastian Manecke for many insightful
discussions and we thank Matthias Beck, Arnau Padrol, and Paco Santos for
helpful remarks on the exposition.  We also thank the two anonymous referees
for their helpful suggestions.
\end{acknowledgment}

\begin{biog}

\item[Kristin Fritsch] 
    is currently finishing her Bachelor degree in mathematics at
    Goethe-Universit\"{a}t Frankfurt. Her mathematical interests are mainly in
    computational finance and she works part time in a forensics department of
    a professional services firm in Frankfurt. Working on this project got her
    in touch with the beautiful interplay of geometry and combinatorics and
    she is now also drawn to discrete mathematics. Kristin 
    loves to travel, in particular when it involves mountains or beautiful cities
    and electronic music in the background.
\begin{affil}
Institut f\"ur Mathematik, Goethe-Universit\"at Frankfurt, Germany\\
kristinfritsch1988@gmail.com
\end{affil}

\item[Janin Heuer] 
    completed her M.S.~in mathematics at the Goethe-Universit\"{a}t Frankfurt
    in March 2019 and is currently pursuing a PhD at Technische Universit\"{a}t
    Braunschweig. Her research interests include polynomial and semidefinite
    optimization. In the future she hopes to be assigned to the Starship
    Enterprise to boldly go where no (wo)man has gone before.
\begin{affil}
    Technische Universit\"at Braunschweig, Institut f\"ur Analysis und
    Algebra, AG Algebra, Universit\"atsplatz 2, 38106 Braunschweig, Germany\\
    janin.heuer@tu-braunschweig.de
\end{affil}

\item[Raman Sanyal] 
    studied at Technische Universit\"{a}t Berlin. After a Miller Research
    Fellowship at UC Berkeley and an assistant professorship at Freie
    Universit\"at Berlin, he joined the mathematics department at the
    Goethe-Universit\"at Frankfurt in 2016. His research is in areas of
    discrete convex geometry and geometric combinatorics, in particular
    polytopes, valuations, and their combinatorics. Raman enjoys
    \{skate,snow,wake\}boarding and reading in a hammock.

\begin{affil}
Institut f\"ur Mathematik, Goethe-Universit\"at Frankfurt, Germany\\
sanyal@math.uni-frankfurt.de
\end{affil}

\item[Nicole Schulz] 
    is currently pursuing a Bachelor degree in mathematics at
    Goethe-Universit\"{a}t Frankfurt.
\begin{affil}
Institut f\"ur Mathematik, Goethe-Universit\"at Frankfurt, Germany\\
nicole.schulz@freenet.de
\end{affil}

\end{biog}
\vfill\eject


\begin{thebibliography}{2}
\bibitem{barvinok}
Barvinok, A. (2002). \textit{A Course in Convexity}.  Graduate Studies in
Mathematics, vol.~54. Providence, RI: American Mathematical Society. 

\bibitem{BR}
Beck, M., Robins, S. (2015). \textit{Computing the Continuous Discretely},
  2nd ed. Undergraduate Texts in Mathematics. New York, NY: Springer.

\bibitem{crt}
Beck, M., Sanyal, R. (2018). \textit{Combinatorial Reciprocity Theorems}.
Graduate Studies in Mathematics, vol.~195. Providence, RI: American
Mathematical Society.

\bibitem{Breuer} 
Breuer, F. (2012). 
Ehrhart $f^\ast$-coefficients of polytopal complexes are non-negative integers.
\textit{Electron. J. Combin.} 19(4): Paper 16, 22 pp. 

\bibitem{BrunsGubeladze}
Bruns, W., Gubeladze, J. (2009). \textit{Polytopes, Rings, and {$K$}-theory}.
Dordrecht, The Netherlands: Springer.

\bibitem{BCKV}
Bump, D., Choi, K., Kurlberg, P., Vaaler, J. (2000).
A local Riemann hypothesis. I. 
\textit{Math. Z.} 233(1):  1--19. 


\bibitem{Triang}
De~Loera, J.~A.,  Rambau, J., Santos, F. (2010). \textit{Triangulations}.
Algorithms and Computation in Mathematics, vol.~25.
Berlin, Germany: Springer.

\bibitem{ehrhartpolynomial} 
Ehrhart, E. (1962). Sur les poly{\`e}dres rationnels homoth{\'e}tiques
{\`a} {$n$}\ dimensions. \textit{C. R. Acad. Sci. Paris}. 254: 616--618.

\bibitem{Gardner}
Gardner, M. (1988). \textit{Hexaflexagons and Other Mathematical Diversions}.
Chicago, IL: Univ.~of Chicago Press.

\bibitem{HKM}
Higashitani, A., Kummer, M., Michalek, M. (2017).
Interlacing Ehrhart polynomials of reflexive polytopes. 
\textit{Selecta Math. (N.S.)} 23(4): 2977--2998. 

\bibitem{KV} 
K\"oppe, M., Verdoolaege, S. (2008).
Computing parametric rational generating functions with a primal Barvinok algorithm. 
\textit{Electron. J. Combin.} 1: Paper 16, 19 pp. 

\bibitem{MS}
Maclagan, D., Sturmfels, B. (2015). \textit{Introduction to Tropical
Geometry}. Graduate Studies in Mathematics, vol.~161. Providence, RI: American
Mathematical Society.

\bibitem{matousek}
Matousek, J. (2002). \textit{Lectures on Discrete Geometry}. 
Graduate Texts in Mathematics, vol.~212. New York, NY: Springer.

\bibitem{StanCyc} 
Stanley, R.~P. (2011). 
Two enumerative results on cycles of permutations.
\textit{European J. Combin.} 32(6):  937--943. 

\bibitem{stanley} 
Stanley, R.~P. (2012). \textit{Enumerative Combinatorics, Volume 1}, 2nd
ed. Cambridge Studies in Advanced Mathematics, vol.~49. New York, NY:
Cambridge Univ. Press.

\bibitem{sullivant} 
Sullivant, S. (2006). Compressed polytopes and statistical disclosure
limitation. \textit{Tohoku Math. J.} 58(3): 433--445. doi.org/10.2748/tmj/1163775139

\bibitem{ziegler}
Ziegler, G.~M. (1995). \textit{Lectures on Polytopes}. Graduate Texts in
Mathematics, vol.~152. New York, NY: Springer.

\end{thebibliography}
\end{document}